\documentclass[12pt]{article}
\usepackage{amssymb}
\usepackage[mathscr]{eucal}
\input{amssym.def}
\input{amssym.tex}

\newtheorem{theorem}{Theorem}
\newtheorem{lemma}{Lemma}
\newtheorem{proposition}{Proposition}
\newtheorem{remark}{Remark}

\newcommand{\GL}{{\rm GL}}

\newcommand{\Hom}{{\rm Hom}}
\newcommand{\ind}{{\rm ind}}
\newcommand{\PGL}{{\rm PGL}}
\newcommand{\proof}{{\it Proof. }}
\newcommand{\sign}{{\rm sign}}

\newcommand{\br}{{\bf r}}
\newcommand{\bs}{{\bf s}}

\def\sq{\hfill$\Box$}

\def\hD{ \mathscr{D} }
\def\hF{ \mathscr{F} }

\def\hK{ \mathscr{K} }
\def\hL{ \mathscr{L} }
\def\hM{ \mathscr{M} }

\def\hSch{ \mathscr{Sch} }

\def\path{{\wp}}

\def\scP{{\mathscr P}}

\newcommand{\A}{{\mathbb A}}

\newcommand{\Af}{{\mathbb A}_{f}}
\newcommand{\Afk}{{\mathbb A}_{f}\otimes k}

\newcommand{\C}{{\mathbb C}}
\newcommand{\F}{{\mathbb F}}

\newcommand{\N}{{\mathbb N}}

\newcommand{\Q}{{\mathbb Q}}

\newcommand{\R}{{\mathbb R}}
\newcommand{\Z}{{\mathbb Z}}

\def\gb{{\mathfrak b}}

\def\go{{\mathfrak o}}

\begin{document}

\title{Shintani Cocycles on $\GL_{n}$}
\author{Richard Hill}
\maketitle

\begin{abstract}
\noindent
The aim of this paper is to define an $n-1$-cocycle
 $\sigma$ on $\GL_{n}(\Q)$ with values in a certain space $\hD$
 of distributions on $\Af^{n}\setminus\{0\}$.
Here $\Af$ denotes the ring of finite ad\`{e}les of $\Q$, and the distributions
 take values in the Laurent series $\C((z_{1},\ldots,z_{n}))$.
This cocycle can be used to evaluate special values of Artin $L$-functions on number fields at negative integers.
The construction generalizes that of Solomon \cite{sol}
 in the case $n=2$.
\end{abstract}
\bigskip

\noindent
2000 Mathematics Subject Classification: 11F75, 11F67.

\section{Introduction}

Let $\Af$ denote the ring of finite ad\`{e}les of $\Q$.
Let $\hD$ be the space of distributions on 
 $\Af^{n}\setminus\{0\}$ with values in
 the Laurent series ring $\C((z_{1},\ldots,z_{n}))$.
In other words
$$
	\hD
	=
	\Hom\Big(\hSch(\Af^{n}\setminus\{0\}),\C((z_{1},\ldots,z_{n})) \Big),
$$
 where $\hSch(\Af^{n}\setminus\{0\})$
 denotes the space of Bruhat-Schwartz functions
 $f:\Af^{n}\setminus\{0\}\to\C$.
The aim of this paper is to define an $n-1$-cocycle
 $\sigma$ on $\GL_{n}(\Q)$ with values in $\hD$.
The construction generalizes that of Solomon \cite{sol}
 in the case $n=2$.
A similar but different cocycle was found by and Solomon and Hu \cite{solhu}
 in the case $n=2,3$.
For $n>3$ the cocycle of Solomon and Hu is only defined on a 
Zariski-open subset of $\GL_{n}(\Q)^{n}$.
Following \cite{sol} and \cite{solhu} we shall
 refer to the cocycle $\sigma$ as the \emph{Shintani cocycle}.

Very briefly, to define the Shintani cocycle $\sigma$
we begin by choosing a non-zero vector $v\in\Q^{n}$.
Given $\alpha_{1},\ldots,\alpha_{n}\in \GL_{n}(\Q)$,
we shall define a cone $C$ in $\Q^{n}$.
Roughly speaking, $C$ will be the set of linear combinations
$\sum x_{i}\alpha_{i}v$ with the $x_{i}\in\Q$ positive.
To be more precise, one must also include some of the faces of $C$.
Then the cocycle is given by
$$
	\sigma(\alpha_{1},\ldots,\alpha_{n})(\varphi)
	=
	\sum_{v\in C} \varphi(v)\exp\left(\sum_{i=1}^{n} v_{i}z_{i}\right),
	\qquad
	\varphi
	\in
	\hSch(\Af^{n}\setminus\{0\}).
$$
Solomon and Hu showed how to make sense of the
right hand side of this formula
as an element of $\C((z_{1},\ldots,z_{n}))$;
we shall explain their method in \S3 below.
The difficulty tackled in this paper,
 is to define the cone $C$ correctly
 in the case where the vectors $\alpha_{i}v$
 are not in general position.
These problems are solved in \S4.
In \S5, we describe the case $n=2$ in detail, and show
how $\sigma$ is related to Solomon's cocycle.

\section{Special values of $L$-functions}

Motivation for studying Shintani-cocycles lies in their relation with
 special values of $L$-functions.
We shall spend a few moments discussing this connection.
Since these matters are the subject of other papers
 (\cite{hida,shintani,sol0,sol,sol2}), the
discussion here is deliberately informal, 
and is independent of the rest of the paper.

Given a totally real number field $k$ of degree $n$,
 one can obtain formulae for the values of abelian $L$-functions
 of $k$ at negative integers by substituting
 units of $k$ (regarded as $n\times n$ integer matrices)
 into the $n-1$-cocycle on $\GL_{n}(\Q)$.
Similar formulae have been obtained
 by Sczech \cite{sczech0,sczech} and Stevens \cite{stevens}.

\paragraph{1. $L$-functions of $\Q$.}
As a first example, we consider the case $n=1$.
In this case we have a (homogeneous) $0$-cocycle $\sigma$ on $\GL_{1}(\Q)$.
Evaluating $\sigma$ at the identity element, we obtain the following
 distribution
$$
	\sigma(1)(\varphi)
	=
	\sum_{v\in \Q,\; v>0}
	\varphi(v)
	\exp(zv)
	\in
	\C((z)),\quad
	\varphi\in \hSch(\Af\setminus\{0\}).
$$
If we substitute for $\phi$ a Dirichlet
character $\chi: \hat\Z \to\C$, extended to be
zero on $\Af\setminus\hat\Z$, then this gives
$$
	\sigma(1)(\chi)
	=
	\sum_{n=1}^{\infty}
	\chi(n)
	\exp(nz).
$$
Differentiating with respect to $z$, we obtain formally (following Euler):
$$
	L(\chi,1-r)
	=
	\left.
	\left(\frac{\partial}{\partial z}\right)^{r-1}
	\Big(\sigma(1)(\chi)\Big)
	\right|_{z=0},
	\qquad
	r\in\N.
$$
To make sense of this equation, let $f\in\N$ be a conductor of $\chi$.
We can group the terms
 into finitely many geometric progressions as follows:
\begin{eqnarray*}
	\sigma(1)(\chi)
	&=&
	\sum_{n=1}^{f}
	\chi(n)
	\exp(nx)
	\bigg(1+\exp(fz)+\exp(2fz)+\cdots \bigg)\\
	&=&
	\sum_{n=1}^{f}
	\chi(n)
	\frac{\exp(nz)}{1-\exp(fz)}.
\end{eqnarray*}
The ratio of exponentials can be expanded in terms of Bernoulli polynomials
 $B_{m}$ as follows.
\begin{eqnarray*}
	\sigma(1)(\chi)
	&=&
	-\sum_{m=0}^{\infty}
	\left(
	\sum_{n=1}^{f}
	\chi(n) B_{m+1}\left(\frac{n}{f}\right)
	\right)
	\frac{(fz)^{m}}{(m+1)!}.
\end{eqnarray*}
This gives the usual expression for $L(\chi,1-r)$ (see for example \S2.3 of \cite{hida}):
\begin{eqnarray*}
	L(\chi,1-r)
	&=&
	-\frac{f^{r-1}}{r}
	\sum_{n=1}^{f}
	\chi(n) B_{r}\left(\frac{n}{f}\right).
\end{eqnarray*}

\paragraph{2. $L$-functions of real quadratic fields.}
Now suppose $k$ is a real quadratic field with
ring of integers $\go$.
For simplicity, we shall assume that $k$ has narrow class number $1$,
 i.e. every non-zero ideal of $\go$ has a totally positive generator.
Indeed, if $u$ denotes a generator for the group of totally positive
 units in $\go$, then each ideal has a unique generator in the following cone:
$$
	C
	=
	\{x+yu: x,y\in \Q,\; x\ge 0,\; y>0\}.
$$
We may therefore express the abelian $L$-functions of $k$ as follows:
$$
	L(\chi,s)
	=
	\sum_{a\in C\cap \go}
	\chi(a)
	N(a)^{-s}.
$$
By choosing an integral basis $\{b_{1},b_{2}\}$, we can regard $k^{\times}$
as a subgroup of $\GL_{2}(\Q)$ and $\go^{\times}$ as a
subgroup of $\GL_{2}(\Z)$.
The special values of $L(\chi,s)$ are encoded in the restriction
of $\sigma$ to $\go^{\times}$.
We shall also identify $\Af^{2}$ with the ring $\Afk$ of finite adeles of $k$.
An element $z\in \Hom_{\Q}(k,\C)$
can be decomposed in terms of the dual basis $\{b_{1}^{*},b_{2}^{*}\}$
as follows:
$$
	z
	=
	z_{1}b_{1}^{*}+z_{2}b_{2}^{*},
	\qquad
	z_{1},z_{2}\in\C.
$$
Evaluating our $1$-cocycle on the fundamental unit $u$,
 we obtain the following the distribution:
$$
	\sigma(1,u)(\varphi)
	=
	\sum_{a\in C}
	\varphi(a)
	\exp(z\cdot a).
$$
Although the right hand side actually converges
 for $z$ in a certain cone in $\Hom_{\Q}(k,\C)$,
 we shall in fact interpret
 it as an element of $\C((z_{1},z_{2}))$ by the method of Solomon and Hu.

Again, substituting a Dirichlet character $\chi$ for $\varphi$,
we obtain:
$$
	\sigma(1,u)(\chi)
	=
	\sum_{a\in C\cap \go}
	\chi(a)
	\exp(z\cdot a).
$$
On the other hand, we can also write $z$ in the form
$$
	z
	=
	t_{1}\tau_{1}+t_{2} \tau_{2},\quad
	t_{1},t_{2}\in\C,
$$
where $\tau_{1}, \tau_{2}:k\to\R$ are the two field embeddings.
With this notation we have
$$
	N(a)^{r}
	=
	\left(
	\frac{\partial}{\partial t_{1}}
	\frac{\partial}{\partial t_{2}}
	\right)^{r}
	\exp(z\cdot a)\Big|_{z=0}.
$$
This gives (formally at least) the following:
$$
	L(\chi,-r)
	=
	\left(
	\frac{\partial}{\partial t_{1}}
	\frac{\partial}{\partial t_{2}}
	\right)^{r}
	\sigma(1,u)(\chi)
	\Big|_{z=0}.
$$
To make sense of this formula,
 we let $f\in\N$ be a (not necessarily minimal) conductor of $\chi$.
We can group the terms of $\sigma(1,u)(\chi)$ into finitely many products of
geometric progressions as before:
\begin{eqnarray*}
	\sigma(1,u)(\chi)
	&=&
	\sum_{a\in \scP\cap \go}
	\chi(a)
	\exp(a\cdot z)
	\sum_{r,s=0}^{\infty} \exp(f(r+su)\cdot z)\\
	&=&
	\sum_{a\in \scP\cap \go}
	\frac{\chi(a)\exp(z\cdot a)}{(1-\exp(z\cdot f))(1-\exp(z\cdot fu))}.
\end{eqnarray*}
Here $\scP$ denotes the half-open parallelogram:
$$
	\scP
	=
	\{x+yu : x,y\in\Q,\; 0\le x <f,\; 0<y\le f\}.
$$
When we expand the ratio of exponentials out
as a power series in $z$, instead of values of
Bernoulli polynomials, the coefficients will instead be
generalized Dedekind sums.
The expansion of such power series in terms of
Dedekind sums and their generalization is
described in \cite{chapman,sol}.
In particular, the reciprocity laws satisfied by these sums
are shown to be consequences of the cocycle relation
satisfied by $\sigma$.
More precisely, we can expand out as follows:
\begin{eqnarray*}
	\sigma(1,u)(\chi)
	&=&
	\sum_{m_{1},m_{2}=0}^{\infty}
	S(m_{1},m_{2},\chi)
	\frac{z_{1}^{m_{1}}}{m_{1}!}
	\frac{z_{2}^{m_{2}}}{m_{2}!}.
\end{eqnarray*}
The coefficients $S(m_{1},m_{2},\chi)$ can be expressed in terms of Dedekind sums.

Using the $r$-th symmetric power of the transition matrix
 from $\{\tau_{1}, \tau_{2}\}$ to $\{b_{1}^{*},b_{2}^{*}\}$,
 we obtain a formula for $L(\chi,-r)$ in terms the numbers
 $S(\chi,m_{1},m_{2})$ with $m_{1}+m_{2}=2r$.

\paragraph{3. Totally Real Fields.}
Let $k$ be a totally real algebraic number field with $[k:\Q]=n$.
Let $H_{\infty}$ denotes the narrow class group of $k$,
 i.e. the group of fractional ideals,
 modulo the principal ideals gererated by totally positive elements.
It was shown by Shintani \cite{shintani}
 (see also \S2.7 of \cite{hida})
 that there are finitely many cones
 $C_{1},\ldots,C_{N}$
 such that any abelian $L$-function of $k$ can be expressed in the form
$$
	L(\chi,s)
	=
	\sum_{[\gb]\in H_{\infty}}
	\sum_{i=1}^{N}
	\sum_{a\in C_{i}\cap \gb^{-1}}
	\chi(a\gb)N(a\gb)^{-s}.
$$
By choosing an integral basis, we can regard $\go^{\times}$
as a subgroup of $\GL_{n}(\Z)$.
The method described above may be used to express the
special values $L(s,-r)$ in terms of the restriction
of $\sigma$ to $\go^{\times}$.

As an example, suppose $k$ has narrow class number 1.
Let $\{u_{1},\ldots,u_{n-1}\}$ be a basis for the group of
totally positive units.
Then we have, for a certain differential operator $\partial$:
$$
	L(\chi,-r)
	=
	\partial^{r}
	\left.\left(
	\sum_{\xi\in S_{n-1}}
	\sign(\xi)\cdot
	\tilde\sigma(u_{\xi(1)},\ldots,u_{\xi(n-1)})
	\right)(\chi)
	\right|_{z=0},
$$
where $\tilde\sigma$ is the corresponding inhomogeneous cocycle:
$$
	\tilde\sigma(\alpha_{1},\ldots,\alpha_{n-1})
	=
	\sigma(1,\alpha_{1},\alpha_{1}\alpha_{2},\ldots,\alpha_{1}\cdots\alpha_{n-1}).
$$

Note that if $k$ has a complex place,
 then its $L$-functions are zero at negative integers.
This can be seen from the functional equation.

\section{Notation and Background Material}

\paragraph{1. The module of cones.}
Let $v_{1},\ldots,v_{r}\in\R^{n}$ be linearly independent vectors.
By the \emph{open cone} of $v_{1},\ldots,v_{r}$, we shall mean
 the set
$$
 C^{o}(v_{1},\ldots,v_{r})
 =
 \left\{
 \sum \lambda_{i} v_{i} :
 \lambda_{1},\ldots,\lambda_{r} > 0
 \right\}.
$$
The \emph{closed cone} $\bar C(v_{1},\ldots,v_{r})$
 is defined similarly but with the inequalities $>$ replaced by $\ge$.
We shall write $\partial C$ for the points of the closed cone which 
 are not in the open cone.
A cone will be called \emph{rational} if
 the vectors $v_{1},\ldots,v_{r}$ are in $\Q^{n}$.
Let $\hK_{\Q}^{o}$ be the abelian group of functions
 $\R^{n}\setminus\{0\} \to\Z$
 generated by the characteristic functions
 of rational open cones.
We shall also write $\hK_{\Q}$ for the group of functions
 $\R^{n} \to\Z$ whose restrictions to $\R^{n}\setminus\{0\}$
 are in $\hK_{\Q}^{o}$.

We shall regard $\hK$ (resp. $\hK_{\Q}$)
 as a left $\GL_{n}(\R)$- (resp. $\GL_{n}(\Q)$-) module
 with the action given by
$$
 (\alpha * c)(v)
 =
 \sign(\det\alpha) \cdot c(\alpha^{-1} v).
$$
The constant functions $\R^{n}\setminus\{0\}\to \Z$ are in $\hK_{\Q}$,
 and form a submodule which we shall denote $\Z(-)$.
The quotient $\hK_{\Q}/\Z(-)$ will be written $\hL_{\Q}$.


\paragraph{2. The Solomon - Hu pairing.}
In \cite{solhu}, Solomon and Hu introduced a
 pairing
$$
 \hL_{\Q}\times \hSch(\A_{f}^{n}\setminus\{0\})\to 
 \C((z_{1},\ldots,z_{n})).
$$
This is defined in several steps.

\underline{Step 1.}
Let $\Z\{\Q^{n}\}$ be the space of all functions $\Q^{n}\to\Z$
 and let $\Z[\Q^{n}]$ be the group ring of the group $\Q^{n}$,
 i.e. the elements of $\Z\{\Q^{n}\}$ of finite support.
One defines a map
 $\Phi:\Z[\Q^{n}]\to \C((z_{1},\ldots,z_{n}))$
 as follows:
$$
 \Phi(A)
 =
 \sum_{w\in\Q^{n}}
 A(w)\exp(w\cdot z).
$$
Here $w\cdot z$ denotes the dot product $w_{1}z_{1}+\ldots+w_{n}z_{n}$.

\underline{Step 2.}
The group ring $\Z[\Q^{n}]$ acts on
 $\Z\{\Q^{n}\}$, $\Z[\Q^{n}]$ and $\C((z_{1},\ldots,z_{n}))$
 and the map $\Phi$ is compatible with these actions.
Define
$$
 \Z\{\Q^{n}\}^{(q)}
 =
 \left\{
 B\in \Z\{\Q^{n}\}:
 \exists A\in\Z[\Q^{n}] \setminus\{0\}
 \hbox{ such that } AB\in \Z[\Q^{n}]
 \right\}.
$$
As $\Z[\Q^{n}]$ is an integral domain,
 it follows that $\Z\{\Q^{n}\}^{(q)}$ is an additive subgroup
 of $\Z\{\Q^{n}\}$.
Furthermore the map
 $\Phi:\Z[\Q^{n}]\to \C((z_{1},\ldots,z_{n}))$
 extends uniquely to $\Z\{\Q^{n}\}^{(q)}$
 in a way which is compatible with the actions of $\Z[\Q^{n}]$.

\underline{Step 3.}
We next define a pairing
 $\hK_{\Q}\times \hSch(\A_{f}^{n})\to \C((z_{1},\ldots,z_{n}))$.
Given $c\in \hK_{\Q}$ and $\varphi\in \hSch(\A_{f}^{n})$,
 we define a function $c\cdot\varphi:\Q^{n}\to\C$ by
$$
	(c\cdot\varphi)(v)
	=
	\cases{
		c(v)\phi(v)	& if $v\ne 0$,
		\smallskip\cr
		0		& if $v=0$.
	}
$$
It turns out that $c\cdot \varphi$ is in $\Z\{\Q^{n}\}^{(q)}$,
and we define
$$
	\langle c,\varphi\rangle
	=
	\Phi(c\cdot\varphi).
$$

We shall describe this more explicitly.
Let $v_{1},\ldots,v_{r}\in\Q^{n}$ be linearly independent
 and let $c$ be the characteristic function of the open cone of 
 $v_{1},\ldots,v_{r}$.
Given $\varphi\in\hSch(\A_{f}^{n})$,
 there is a lattice $L\subset \Q^{n}$ such that
 $\varphi$ is invariant under translation by $L$.
By multiplying the vectors $v_{1},\ldots,v_{n}$ by natural numbers if 
 necessary, we may assume $v_{1},\ldots,v_{n}\in L$.
Let
$$
 \scP
 =
 \{x_{1}v_{1}+\ldots+x_{n}v_{n}:
 x_{1},\ldots,x_{n}\in (0,1]\}.
$$
We have $(1-[v_{1}])\ldots(1-[v_{r}])(c\cdot\varphi)=p\cdot\varphi$,
 where $p$ is the characteristic function of the parallelotope $\scP$.
As $p\cdot\varphi$ has finite support,
 it follows that $c\cdot\phi$ is in $\Z\{\Q^{n}\}^{(q)}$,
 and the pairing is given by:
$$
 \langle c,\varphi\rangle
 =
 \frac{1}{1-\exp(v_{1} \cdot z)}
 \ldots
 \frac{1}{1-\exp(v_{r} \cdot z)}
 \sum_{w\in\scP\cap \Q^{n}}
 \varphi(w) \exp(w\cdot z).
$$

\underline{Step 4.}
Let $c$ be the constant function on $\Q^{n}\setminus\{0\}$
 with value 1 and let $\varphi\in \hSch(\Af^{n})$.
Then for any non-zero vector
 $v\in\Q^{n}$ such that $\varphi$ is periodic modulo $v$,
 we have
$$
 \Big((1-[v])(c\cdot\varphi)\Big)(w)
 =
 \left\{
 \begin{array}{ll}
     \varphi(0) & \hbox{if } w=-v \\
     -\varphi(0) & \hbox{if } w=0 \\
     0 & \hbox{otherwise.}
 \end{array}
 \right.
$$
It follows that the constant functions are orthogonal
to the subspace $\hSch(\Af^{n}\setminus\{0\})$.
Hence the pairing factors through to give
 $\hL_{\Q}\times\hSch(\Af^{n}\setminus\{0\})\to \C((z_{1},\ldots,z_{n}))$.

\paragraph{3. The Result.}
To obtain a cocycle with values in $\hD$,
 it is sufficient, using the pairing defined above,
 to construct one with values in $\hL_{\Q}$.
In the case $n=2$, Solomon's cocycle $s$ (with values in $\hL_{\Q}$)
 is defined as follows.
For $\alpha,\beta\in\GL_{2}(\Q)$,
$$
 s(\alpha,\beta)(v)
 =
 \left\{
 \begin{array}{ll}
     \sign\det (\alpha e_{1},\beta e_{1})
     &
     \hbox{if $v\in C^{o}(\alpha e_{1},\beta e_{1})$,}\\
     \frac{1}{2}
     \sign\det (\alpha e_{1},\beta e_{1})
     &
     \hbox{if $v\in \partial C(\alpha e_{1},\beta e_{1})$,}\\
     0
     &
     \hbox{otherwise.}
 \end{array}
 \right.
$$
Here $e_{1}=\left(\matrix{1\cr 0}\right)$.
One could of course replace $e_{1}$ by any other non-zero
 vector to obtain a cohomologous cocycle.
In the case that $\alpha e_{1},\beta e_{1},\gamma e_{1}$
 are in general position, it is easy to see that $s$ satisfies the 
 cocycle relation:
$$
 s(\beta,\gamma)
 -
 s(\alpha,\gamma)
 +
 s(\alpha,\beta)
 =
 0\;\;
 \hbox{modulo constant functions}.
$$
If $\alpha e_{1},\beta e_{1}$ are not linearly independent
 then the above definition makes no sense and we instead define
 $s(\alpha,\beta)=0$.
With this completed definition the cocycle relation remains true
 modulo the kernel of the Solomon-Hu pairing.

Naively one would expect to generalize the cocycle $s$ above by 
defining for $\alpha_{1},\ldots,\alpha_{n}\in \GL_{n}(\Q)$:
$$
 s(\alpha_{1},\ldots,\alpha_{n})(v)
 =
 \left\{
 \begin{array}{ll}
     \sign\det (\alpha_{1} e_{1},\ldots,\alpha_{n} e_{1})
     &
     \hbox{if $v\in C^{o}(\alpha_{1} e_{1},\ldots,\alpha_{n} e_{1})$,}\\
     0
     &
     \hbox{if $v\notin \bar C(\alpha_{1} e_{1},\ldots,\alpha_{n} e_{1})$,}
 \end{array}
 \right.
$$
Indeed as long as $v,\alpha_{1}e_{1},\ldots,\alpha_{n}e_{1}$ are in 
 general position, the above definition makes sense and a similar
 cocycle relation is satisfied.
The difficulties are (a) how to define $s$
 when $\alpha_{1}e_{1},\ldots,\alpha_{n}e_{1}$
 are linearly dependent, and
 (b) how to define $s$ when
 $v\in \partial C(\alpha_{1} e_{1},\ldots,\alpha_{n} e_{1})$
 without losing the cocycle relation.
Both these problems are solved by the same method in this paper.

\section{Definition of the Shintani Cocycle}

\paragraph{1. A relation between signs of determinants}

By an \emph{ordered field} we shall mean a (commutative) field
 $\F$ equipped with a total ordering $>$ satisfying the condition:
\begin{itemize}
    \item
    $\forall x,y,z\in \F$ if $x>y$ then $x+z>y+z$;
    \item
    $\forall x,y,z\in \F$ if $x>y$ and $z>0$ then $xz>yz$.
\end{itemize}
Fix an ordered field $\F$ and define for $x\in \F^{\times}$:
$$
 \sign(x)
 =
 \left\{
 \begin{array}{ll}
     1& x>0,\\
     -1 & x<0.
 \end{array}
 \right.
$$
A set of vectors in $\F^{n}$ will be said to be
 \emph{in general position} if every subset with no more than $n$ elements
 is linearly independent.
Given vectors $v_{0},\ldots,v_{n}\in \F^{n}$ in general position,
 there are non-zero scalars $\lambda_{0},\ldots,\lambda_{n}\in \F$
 such that $\lambda_{0}v_{0}+\ldots+\lambda_{n}v_{n}=0$.
We define
$$
 d(v_{0},\ldots,v_{n})
 =
 \left\{
 \begin{array}{ll}
     (-1)^{i}\sign\det(v_{1},\ldots,\hat v_{i},\ldots,v_{n}) &
     \hbox{if $\lambda_{0},\ldots,\lambda_{n}$ all have the same 
     sign,}\\
     0 & \hbox{otherwise.}
 \end{array}
 \right.
$$

\begin{proposition}
    \label{dprops}
    \begin{itemize}
	\item[(i)]
	$d(v_{0},\ldots,v_{n})$ is well defined (i.e. independent of $i$).
	\item[(ii)]
	For any permutation $\xi$ we have
	$d(v_{\xi(0)},\ldots,v_{\xi(n)})
	 =\sign(\xi)d(v_{0},\ldots,v_{n})$.
	\item[(iii)]
	If $\lambda_{0},\ldots,\lambda_{n}\in \F^{>0}$
	 then $d(\lambda_{0}v_{0},\ldots,\lambda_{n}v_{n})
	 = d(v_{0},\ldots,v_{n})$.
	\item[(iv)]
	For any $\alpha\in\GL_{n}(\F)$ we have
	$d(\alpha v_{0},\ldots,\alpha v_{n}) =
	 \sign(\det\alpha) \cdot d(v_{0},\ldots,v_{n})$.
	\item[(v)]
	Let $\F'$ be another ordered field and let
	 $\iota:\F\hookrightarrow \F'$ be an order-preserving field homomorphism.
	Then for $v_{0},\ldots,v_{n}$ in general position in $\F^{n}$
	 we have $d(\iota v_{0},\ldots,\iota v_{n}) =
	 d(v_{0},\ldots,v_{n})$.	
	\item[(vi)]
	If $v_{1},\ldots,v_{n+2}\in \F^{n}$ are in general position then
	$$
	 \sum_{i=0}^{n+2} (-1)^{i}d(v_{1},\ldots,\hat v_{i},\ldots,v_{n+2})
	 =
	 0.
	$$
    \end{itemize}
\end{proposition}

\proof
(i)
The case that $d(v_{0},\ldots,v_{n})=0$ is clearly well defined,
 so assume $v_{0}=-\sum_{j=1}^{n}\lambda_{j} v_{j}$
 with $\lambda_{j}>0$.
We have by elementary properties of determinants:
\begin{eqnarray*}
    (-1)^{i}\sign\det(v_{0},v_{1},\ldots,\hat v_{i},\ldots,v_{n})
    &=&
    (-1)^{i}
    \sign\det(-\lambda_{i}v_{i},v_{1},\ldots,\hat v_{i},\ldots,v_{n})\\
    &=&
    (-1)^{i-1}
    \sign\det(v_{i},v_{1},\ldots,\hat v_{i},\ldots,v_{n})\\
    &=&
    \sign\det(v_{1},\ldots,v_{n}).
\end{eqnarray*}

(ii)
This follows from (i) in the case that $\xi$ is an adjacent
 transposition $(i \quad i+1)$.
The general case follows since the adjacent transpositions generate
 the group of all permutations.

Parts (iii), (iv) and (v) follow immediately from the definition.

(vi)
By (ii), (iii) and (iv) we may reduce the general case
 to the case  that $v_{1},\ldots,v_{n}$
 are the standard basis elements $e_{1},\ldots,e_{n}$ in $\F^{n}$
 and $v_{n+1}$ and $v_{n+2}$ are of the form
$$
 v_{n+1}
 =
 \left(
 \matrix{
 \left.\matrix{1\cr \vdots \cr 1}\right\}r\cr
 \matrix{-1\cr \vdots \cr -1}\quad}
 \right);
 \quad
 v_{n+2}
 =
 \left(\matrix{x_{1}\cr \vdots \cr x_{r} \cr x_{r+1} \cr \vdots \cr 
 x_{n}}\right),
$$
 with $x_{1}>\ldots>x_{r}$ and $x_{r+1}<\ldots<x_{n}$.
For simplicity we assume $0<r<n$;
 the cases $r=0$ and $r=n$ may be handled similarly.

Case 1. Assume $i\le r$.
If we have
$$
 v_{n+2}
 =
 \sum_{j=1}^{n+1} \lambda_{j}v_{j},
$$
 then this implies
$$
 \lambda_{j}
 =
 \left\{
 \begin{array}{ll}
     x_{j}-x_{i} & \hbox{if $j\le r$, $j\ne i$,}\\
     x_{j}+x_{i} & \hbox{if $r<j\le n$.}\\
     x_{i} & \hbox{if $j=n+1$.}
 \end{array}
 \right. 
$$
For all the coefficients $\lambda_{j}$ to be negative
 we require $i=1$, $x_{1}<0$ and $x_{1}+x_{n}<0$.
From this we may deduce that
$$
 \sum_{i=1}^{r}(-1)^{i}d(v_{1},\ldots,\hat v_{i},\ldots,v_{n})
 =
 \left\{
 \begin{array}{ll}
     1 & \hbox{if $x_{1}<0$ and $x_{1}+x_{n}<0$,}\\
     0 & \hbox{otherwise.}
 \end{array}
 \right.
$$
A similar calculation shows that
$$
 \sum_{i=r+1}^{n}
 (-1)^{i} d(v_{1},\ldots,\hat v_{i},\ldots,v_{n})
 =
 \left\{\begin{array}{ll}
     -1 & \hbox{if $x_{n}>0$ and $x_{1}+x_{n}<0$,}\\
     0 & \hbox{otherwise.}
 \end{array}
 \right.
$$
Finally we have $d(v_{1},\ldots,v_{n+1})=0$ and
$$
 (-1)^{n+1} d(v_{1},\ldots, v_{n},\ldots,v_{n+2})
 =
 \left\{\begin{array}{ll}
     -1 & \hbox{if $x_{1}<0$ and $x_{n}<0$,}\\
     0 & \hbox{otherwise.}
 \end{array}
 \right.
$$
Adding everything up we obtain the result.
\sq
\medskip

\paragraph{2. A relation between cone functions}
Given a basis $\{v_{1},\ldots,v_{n}\}$ of $\F^{n}$
 we define a function $c(v_{1},\ldots,v_{n}):\F^{n}\to \Z$ by
$$
 c(v_{1},\ldots,v_{n})(w)
 =
 \left\{
 \begin{array}{ll}
     \sign\det(v_{1},\ldots,v_{n}) &
     \hbox{if }v_{1}^{*}(w),\ldots,v_{n}^{*}(w)>0,\\
     0 & \hbox{otherwise.}
 \end{array}
 \right.
$$
Thus up to a sign, $c(v_{1},\ldots,v_{n})$ is the characteristic 
 function of the open cone of $v_{1},\ldots,v_{n}$.
This may be expressed in terms of the function $d$:
\begin{equation}
   \label{cd}
   c(v_{1},\ldots,v_{n})(w)
   =
   (-1)^{n}d(v_{1},\ldots,v_{n},-w).
\end{equation}
The functions $c$ and $d$ satisfy the following cocycle relation.

\begin{proposition}
    \label{cprop}
    If $v_{0},\ldots,v_{n},w\in \F^{n}$ are in general position then
    $$
     \sum_{i=0}^{n} 
     (-1)^{i}c(v_{0},\ldots,\hat{v_{i}},\ldots,v_{n})(w)
     =
     d(v_{0},\ldots,v_{n}).
    $$
\end{proposition}

\proof
This follows immediately from Proposition \ref{dprops} (vi)
 and (\ref{cd}).
\sq
\medskip

\paragraph{3.}
Consider the local field $\F=\R((\epsilon_{1}))\ldots((\epsilon_{n}))$.
Every element of $\F$ may be expressed as
$$
 \sum_{\br\in\Z^{n}} a_{\br}\epsilon^{\br},
$$
 where $\epsilon^{\br}=\epsilon_{1}^{r_{1}}\cdots \epsilon_{n}^{r_{n}}$.
The coefficients $a_{\br}$ are in $\R$.
We shall order the multi-indices $\br\in\Z^{n}$ lexicographically,
 so $\br<\bs$ if and only if there is an $i\in\{1,2,\ldots,n\}$ such that
$$
 r_{i}<s_{i}
 \hbox{ and }
 \forall j>i, \; r_{j}=s_{j}.
$$
Using this ordering we may define the leading term of a non-zero 
 element of $\F$ to be the non-zero monomial $a_{\br}\epsilon^{\br}$ for 
 which $\br$ is smallest.
An element of $\F^{\times}$ will be said to be positive (resp. negative)
 if its leading term $a_{\br}\epsilon^{\br}$ has positive (resp. negative)
 coefficient.
For $f,g\in \F$ we define $f>g$ if and only if $f-g$ is positive.
Under this ordering $\epsilon_{1}$ is positive
 but smaller than every positive real number.
For $i=1,\ldots,n-1$, the element $\epsilon_{i+1}$ is
 positive but smaller than every power of $\epsilon_{i}$.

We shall also use the field
 $\F'=\R((\epsilon_{0}))\ldots((\epsilon_{n}))$,
 ordered in an analogous way.
We have $n+1$ order-preserving field embeddings
 $\iota_{i}:\F\hookrightarrow \F'$ defined by
$$
 (\iota_{i}f)(\epsilon_{0},\ldots,\epsilon_{n})
 =
 f(\epsilon_{0},\ldots,\hat\epsilon_{i},\ldots,\epsilon_{n}).
$$

\paragraph{4.}
Define for $i=1,\ldots,n$:
$$
 b(\epsilon_{i})
 =
 \left(\matrix{1\cr\epsilon_{i}\cr\vdots\cr\epsilon_{i}^{n-1}}\right).
$$
We shall regard $b(\epsilon_{i})$ as an element of $\F^{n}$.

\begin{lemma}
    \label{general}
    For any $\alpha_{0},\ldots,\alpha_{n}\in \GL_{n}(\R)$
     and any $w\in \R^{n}\setminus\{0\}$,
     the set
     $\{\alpha_{0}b(\epsilon_{0}),\ldots,\alpha_{n}b(\epsilon_{n}),w\}$
     is in general position in $\F'{}^{n}$.
\end{lemma}

\proof
Regarding $b$ as a function $\R\to\R^{n}$, we note that the values of
 $\alpha_{i}b(\epsilon_{i})$ span $\R^{n}$.
We may therefore choose $\epsilon_{1},\ldots,\epsilon_{n}\in\R$
 so that $\{\alpha_{1}b(\epsilon_{1}),\ldots,\alpha_{n}b(\epsilon_{n})\}$
 is a basis of $\R^{n}$.
Hence $\det(\alpha_{1}b(\epsilon_{1}),\ldots,\alpha_{n}b(\epsilon_{n}))$
 is a non-zero function of $\epsilon_{1},\ldots,\epsilon_{n}$,
 so is a non-zero element of $\F$.
It follows that $\{\alpha_{1}b(\epsilon_{1}),\ldots,\alpha_{n}b(\epsilon_{n})\}$
 is a basis of $\F^{n}$.
A similar argument shows that for any $j$,
 $\{w,\alpha_{i}b(\epsilon_{i}):i\ne j\}$
 is also a basis of $\F^{n}$.
\sq
\medskip

\paragraph{5.}
We now define our cocycle.
Let $\hF$ denote the space of all functions 
 $\varphi:\R^{n}\setminus\{0\}\to \Z$.
We let $\GL_{n}(\R)$ act on $\hF$ by:
$$
 (\alpha * \varphi)(w)
 =
 \sign\det\alpha \cdot \varphi(\alpha^{-1}w),\quad
 \alpha\in\GL_{n}(\R),\;
 \varphi\in \hF,\;
 w\in\R^{n}\setminus\{0\}.
$$
The constant functions in $\hF$ form a submodule, which we shall 
 denote $\Z(-)$.
We shall write $\hM$ for the quotient.
We shall describe a cocycle $\sigma\in H^{n-1}(\GL_{n}(\R),\hM)$.

For $\alpha_{1},\ldots,\alpha_{n}\in\GL_{n}(\R)$
 and $w\in\R^{n}\setminus\{0\}$,
 we define
$$
 \sigma(\alpha_{1},\ldots,\alpha_{n})(w)
 =
 c(\alpha_{1}b(\epsilon_{1}),\ldots,\alpha_{n}b(\epsilon_{n}))(w).
$$

\begin{proposition}
    \label{cocycle}
    (i)
    For $\alpha_{0},\ldots,\alpha_{n}\in\GL_{n}(\R)$
     and $w\in\R^{n}\setminus\{0\}$ we have
    $$
     \sum_{i=0}^{n} 
     (-1)^{i}\sigma(\alpha_{0},\ldots,\hat\alpha_{i},\ldots,\alpha_{n})(w)
     =
     d(\alpha_{0}b(\epsilon_{0}),\ldots,\alpha_{n}b(\epsilon_{n})).
    $$
    (ii)
    For $\beta,\alpha_{1},\ldots,\alpha_{n}\in\GL_{n}(\R)$ we have
    $$
     \sigma(\beta\alpha_{1},\ldots,\beta\alpha_{n})
     =
     \beta * \sigma(\alpha_{1},\ldots,\alpha_{n}).
    $$
\end{proposition}

\proof
(i)
We have by definition
$$
 \sigma(\alpha_{0},\ldots,\hat\alpha_{i},\ldots,\alpha_{n})(w)
 =
 c(\alpha_{0}b(\epsilon_{1}),\ldots,
 \alpha_{i-1}b(\epsilon_{i}),\alpha_{i+1}b(\epsilon_{i+1}),
 \ldots,\alpha_{n}b(\epsilon_{n}))(w).
$$ 
Applying the the order-preserving map $\iota_{i}:\F\to \F'$
 we have
$$
 \sigma(\alpha_{0},\ldots,\hat\alpha_{i},\ldots,\alpha_{n})(w)
 =
 c(\alpha_{0}b(\epsilon_{0}),\ldots,
 \alpha_{i-1}b(\epsilon_{i-1}),\alpha_{i+1}b(\epsilon_{i+1}),
 \ldots,\alpha_{n}b(\epsilon_{n}))(w).
$$ 
The result now follows from Lemma \ref{general}
 and Proposition \ref{cprop}.

(ii)
This follows from Propsition \ref{dprops} (iv) and (\ref{cd}).
\sq
\medskip

The proposition shows that $\sigma$ represents an element of
 $H^{n-1}(\GL_{n}(\R),\hM)$.
The short exact sequence
$$
 0 \to \Z(-)\to \hF \to \hM \to 0,
$$
 gives rise to a connecting homomorphism
 $\partial:H^{n-1}(\GL_{n},\hM)\to H^{n}(\GL_{n},\Z(-))$.
The proposition also shows that $\partial\sigma$
 is given by the $n$-cocycle
$$
 \tau(\alpha_{0},\ldots,\alpha_{n})
 =
 d(\alpha_{0}b(\epsilon_{0}),\ldots,\alpha_{n}b(\epsilon_{n})).
$$

\paragraph{Aside.}
In this context it is worth recording the following long exact sequence:
$$
 \ldots
 \to
 H^{r}(\GL_{n}(\R),\Z(-))
 \to
 H^{r}(\GL_{n-1}(\R),\Z(-))
 \to
 H^{r}(\GL_{n}(\R),\hM)
 \to
 H^{r+1}(\GL_{n}(\R),\Z(-))
 \to
 \ldots.
$$

\proof
We need only show that $H^{r}(\GL_{n-1}(\R),\Z(-))$
 is canonically isomorphic to $H^{r}(\GL_{n}(\R),\hF)$.
Consider the \emph{mirabolic} subgroup:
$$
 P
 =
 \{\alpha\in\GL_{n}(\R) : \alpha e_{1}=e_{1}\}.
$$
We have $\hF=\ind_{P}^{\GL_{n}(\R)} \Z(-)$.
Hence by Shapiro's Lemma,
$$
 H^{r}(\GL_{n}(\R),\hF)
 =
 H^{r}(P,\Z(-)).
$$
The group extension
$$
 1 \to \R^{n-1}\to P \to \GL_{n-1}(\R) \to 1,
$$
gives rise to the spectral sequence
$$
 H^{p}(\GL_{n-1}(\R),H^{q}(\R^{n-1},\Z(-))) \Rightarrow H^{p+q}(P,\Z(-)).
$$
The result now follows since
$$
 H^{q}(\R^{n-1},\Z(-))
 =
 \left\{
 \begin{array}{lll}
     \Z(-)& \hbox{if} & q=0,\\
     0 & &
     q>0.
 \end{array}
 \right.
$$
\sq
\medskip

\paragraph{6.}
We finally show that the values of $\sigma$ are actually in the 
 module of cones.


\begin{theorem}
    (i) For $\alpha_{1},\ldots,\alpha_{n}\in \GL_{n}(\R)$
     we have $\sigma(\alpha_{1},\ldots,\alpha_{n})\in\hL_{\R}$.\\
    (ii) For $\alpha_{1},\ldots,\alpha_{n}\in \GL_{n}(\Q)$
     we have $\sigma(\alpha_{1},\ldots,\alpha_{n})\in\hL_{\Q}$.
\end{theorem}

\proof
Note that $\hK_{\R}$ is closed under pointwise multiplication
 of functions.
Hence to prove the first part of the proposition,
 it is sufficient to show that, for any linear form
 $\phi:\F^{n}\to \F$, the set
$$
 S=\{w\in\R^{n} : \phi(w) >0\}
$$
 is a finite disjoint union of open cones.
The restriction $\phi:\R^{n}\to \F$ is $\R$-linear.
We may write $\phi$ as
$$
 \phi
 =
 \sum_{\br} \phi_{\br}\epsilon^{\br},
$$
 with $\phi_{\br}:\R^{n}\to \R$ linear forms. 
In this sum $\br$ runs over the multipowers of the $\epsilon_{i}$.
We may therefore decompose $S$ into disjoint subsets:
$$
 S
 =
 \bigcup_{\br}S_{\br},
$$
 where
$$
 S_{\br}
 =
 \{w\in \R^{n}:
 \phi_{\br}(w)>0
 \hbox{ and for all $\bs<\br$, }
 \phi_{\bs}(w)=0
 \}
$$
Each non-empty $S_{\br}$ is an open half-subspace,
 and is hence a finite disjoint union of cones.
It remains to show that only finitely many $S_{\br}$ are non-empty.
If $S_{\br}$ and $S_{\bs}$ are both non-empty and $\bs<\br$
 then $S_{\br}$ is contained in the boundary of the closure of $S_{\bs}$
 and is therefore of strictly smaller dimension.

This proves the first part of the proposition.
Now assume $\alpha_{1},\ldots,\alpha_{n}\in\GL_{n}(\Q)$.
It follows that the basis vectors $\alpha_{i}b(\epsilon_{i})$
 are in $\Q(\epsilon)$.
From this it follows that $\phi_{\br}:\Q^{n}\to \Q$.
Hence the sets $S_{\br}$ may be decomposed into rational cones.
\sq
\medskip

\begin{remark}
One could define a $K$-cone for any subfield $K$ of $\R$
 and obtain a generalization of the above proposition.
\end{remark}

\section{Comparison with previous results}

We shall describe $\sigma$ in the case $n=2$ and then give
 a coboundary relating it to Solomon's cocycle $s$.

\paragraph{1.}
As $\sigma$ is homogeneous, we need only calculate $\sigma(1,\alpha)$
 for $\alpha=\left(\matrix{a&b\cr c&d}\right) \in \GL_{2}(\R)$.
Let $w=\left(\matrix{x\cr y}\right)\in\R^{2}\setminus\{0\}$.
Recall that to calculate $\sigma(1,\alpha)(w)$,
 we express $w$ in the form
$$
 w
 =
 x'\left(\matrix{1\cr\epsilon_{1}}\right)
 +
 y'\alpha\left(\matrix{1\cr\epsilon_{2}}\right),
 \quad
 x',y'\in\R((\epsilon_{1}))((\epsilon_{2})).
$$
To simplify notation consider the matrix
$$
 M
 =
 \left(
 \left(\matrix{1\cr\epsilon_{1}}\right),
 \alpha\left(\matrix{1\cr\epsilon_{2}}\right)
 \right)
 =
 \left(
 \matrix{1 & a+b\epsilon_{2}\cr
 \epsilon_{1}& c+d\epsilon_{2}}
 \right) .
$$
We have
$$
 \left(\matrix{x'\cr y'}\right)
 =
 M^{-1}\left(\matrix{x\cr y}\right).
$$
The cocycle is given by the formula:
$$
 \sigma(1,\alpha)\left(\matrix{x\cr y}\right)
 =
 \left\{
 \begin{array}{ll}
     \sign(\det M)
     & \hbox{if $x'$ and $y'$ are both positive in 
                  $\R((\epsilon_{1}))((\epsilon_{2}))$,}\cr
     0
     & \hbox{otherwise.}
 \end{array}
 \right.
$$
After solving these inequalities we obtain:

\begin{proposition}
    \label{example}
    (i) Let $\alpha=\left(\matrix{a&b\cr 0&c}\right)$.
    \begin{itemize}
	\item
	If $a>0$ and $c>0$ then $\sigma(1,\alpha)=0$.
	\item
	If $a>0$ and $c<0$ then
	$
	 \sigma(1,\alpha)\left(\matrix{x\cr y}\right)
	 =
	 \left\{
	 \begin{array}{ll}
	     -1 & \hbox{if $x>0$ and $y=0$,}\cr
	     0 & \hbox{otherwise.}
	 \end{array}
	 \right.
	$
	\item
	If $a<0$ and $c>0$ then
	$
	 \sigma(1,\alpha)\left(\matrix{x\cr y}\right)
	 =
	 \left\{
	 \begin{array}{ll}
	     1 & \hbox{if $y>0$,}\cr
	     0 & \hbox{otherwise.}
	 \end{array}
	 \right.
	$
	\item
	If $a<0$ and $c<0$ then
	$
	 \sigma(1,\alpha)\left(\matrix{x\cr y}\right)
	 =
	 \left\{
	 \begin{array}{ll}
	     1 & \hbox{if $y>0$ or if $y=0$ and $x<0$}\cr
	     0 & \hbox{otherwise.}
	 \end{array}
	 \right.
	$
    \end{itemize}
    (ii)
    Let
     $\alpha
     =
     \left(\matrix{a&b\cr 0&c}\right)
     \left(\matrix{0&1\cr 1&0}\right)
     \left(\matrix{1&d\cr 0&1}\right)$.
    \begin{itemize}
	\item
	If $a>0$ and $c>0$ then
	$
	 \sigma(1,\alpha)\left(\matrix{x\cr y}\right)
	 =
	 \left\{
	 \begin{array}{ll}
	     1 & \hbox{if $y>0$ and $cx-by>0$,}\cr
	     0 & \hbox{otherwise.}
	 \end{array}
	 \right.
	$
	\item
	If $a>0$ and $c<0$ then
	$
	 \sigma(1,\alpha)\left(\matrix{x\cr y}\right)
	 =
	 \left\{
	 \begin{array}{ll}
	     -1 & \hbox{if $y\le 0$ and $cx-by<0$,}\cr
	     0 & \hbox{otherwise.}
	 \end{array}
	 \right.
	$
	\item
	If $a<0$ and $c>0$ then
	$
	 \sigma(1,\alpha)\left(\matrix{x\cr y}\right)
	 =
	 \left\{
	 \begin{array}{ll}
	     1 & \hbox{if $y> 0$ and $cx-by\le 0$,}\cr
	     0 & \hbox{otherwise.}
	 \end{array}
	 \right.
	$
	\item
	If $a<0$ and $c<0$ then
	$
	 \sigma(1,\alpha)\left(\matrix{x\cr y}\right)
	 =
	 \left\{
	 \begin{array}{ll}
	     -1 & \hbox{if $y\le 0$ and $cx-by\le 0$,}\cr
	     0 & \hbox{otherwise.}
	 \end{array}
	 \right.
	$
    \end{itemize}
\end{proposition}

\proof
To give an impression of how to do this calculation,
 we shall prove (i) in the case $a,c<0$.
The other cases are left to the reader.
We have
 $M=\left(\matrix{1&a+b\epsilon_{2}\cr \epsilon_{1}& 
 c\epsilon_{2}}\right)$.
Hence $\det M=c\epsilon_{2}-a\epsilon_{1}-b\epsilon_{1}\epsilon_{2}$.
The leading term of $\det M$ is $-a\epsilon_{1}$, which is positive.
Therefore $\det M >0$.
Furthermore
$$
 x'
 =
 \frac{1}{\det M}
 \left(-ay + (cx-by)\epsilon_{2}\right),
 \quad
 y'
 =
 \frac{1}{\det M}
 \left(-\epsilon_{1} x +y\right).
$$
For $y'$ to be positive we require either $y>0$ or $y=0$ and $x<0$.
In both of these cases $x'$ is also positive.
\sq
\medskip

\paragraph{2.}
In \cite{sol} Solomon obtained a cocycle on $\PGL_{2}$
 rather than on $\GL_{2}$; however the values of the cocycle
 were in $\frac{1}{2}\hL$ rather than in $\hL$.
This cocycle $s\in Z^{1}(\PGL_{2}(\R),\frac{1}{2}\hL)$
 is defined as follows:
$$
 s(\alpha,\beta)(w)
 =
 \left\{
 \begin{array}{ll}
     \sign\det(\alpha e_{1},\beta e_{1}) &
     \hbox{if $\{\alpha e_{1},\beta e_{1}\}$ is a basis of $\R^{2}$
     and $w\in C^{o}(\alpha e_{1},\beta e_{1})$},\\
     \frac{1}{2}\sign\det(\alpha e_{1},\beta e_{1}) &
     \hbox{if $\{\alpha e_{1},\beta e_{1}\}$ is a basis of $\R^{2}$
     and $w\in \partial C(\alpha e_{1},\beta e_{1})$},\\
     0 & \hbox{otherwise.}
 \end{array}
 \right.
$$
This is related to $\sigma$ by the coboundary:
$$
 (\sigma-s)(\alpha,\beta)
 =
 \alpha*\tau-\beta*\tau,
$$
 where
$$
 \tau\left(\matrix{x\cr y}\right)
 =
 \left\{
 \begin{array}{ll}
     \frac{1}{2} &
     \hbox{if $y=0$ and $x>0$,}\\
     0 &
     \hbox{otherwise.}
 \end{array}
 \right.
$$

\bigskip
\bigskip

\noindent
{\tt rih@math.ucl.ac.uk}\\

\noindent
Department of Mathematics,\\
University College London,\\
Gower Street,\\
London WC1E 6BT\\
UK.


\begin{thebibliography}{10}

	\bibitem{chapman}
	R. Chapman,
	``Reciprocity Laws for generalized higher dimensional Dedekind Sums''
	\emph{Acta Arithmetica} XCIII.2 (2000) 189-199.
	
	\bibitem{hida}
	H. Hida,
	``Elementary theory of $L$-functions and Eisenstein series''
	\emph{London Mathematical Society Student Texts} 26,
	Cambridge University Press 1993.
	
	\bibitem{solhu}
	S. Hu and D. Solomon,
	``Properties of higher-dimensional Shintani 
	generating functions and cocycles on ${\rm PGL}_3(\Q)$'',
	\emph{Proc. L.M.S.} 82 no. 1 (2001).
	
	\bibitem{sczech0}
	R. Sczech, 
	``Eisenstein cocycles for $\GL_{2}(\Q)$ and values of $L$-functions 
	in real quadratic fields'' 
	\emph{Comment. Math. Helvetici} 67 (1992) 363-382 .

	\bibitem{sczech}
	R. Sczech,
	``Eisenstein group cocycles for $\GL_{n}$ and values of $L$-functions''
	\emph{Invent. Math.} 113(3) 581-616, 1993.
	
	\bibitem{shintani}
	T. Shintani,
	``On evaluation of zeta functions of totally real algebraic number fields
	at non-positive integers''
	\emph{J. Fac. Sci. Univ. Tokyo Sect. IA Math.} 23 (1976), no. 2, 393-417.
	
	\bibitem{sol0}
	D. Solomon,
	``$p$-adic limits of Shintani generating functions for a real quadratic field.''
	\emph{J. Number Theory} 59 (1996), no. 1, 119--158.

	\bibitem{sol}
	D. Solomon,
	``Algebraic Properties of Shintani Generating Functions:
	Dedekind Sums and Cocycles on PGL(2,$\Q$)'',
	\emph{Compositio Math.} 12 no. 3 (1998).
	
	\bibitem{sol2}
	D. Solomon,
	``The Shintani cocycle. II. Partial $\zeta$-functions,
	cohomologous cocycles and $p$-adic interpolation''
	\emph{J. Number Theory} 75 (1999), no. 1, 53-108.
		
	\bibitem{stevens}
	G. Stevens,
	``The Eisenstein measure and real quadratic fields''
	in \emph{Th\'eorie des nombres (Quebec, PQ, 1987)}
	pages 887-927 de Gruyter, Berlin, 1989.


\end{thebibliography}
\end{document}